\documentclass[12pt]
{amsart}

\usepackage[english]{babel}
\usepackage{amsmath,amssymb}%{amscd}
\usepackage[matrix,arrow]{xy}
\usepackage[active]{srcltx} %SRC Specials for DVI Searching   

\makeatletter
\@namedef{subjclassname@2010}{%
  \textup{2010} Mathematics Subject Classification}
\makeatother

\newtheorem{thm}{Theorem}[section]
\newtheorem{lem}[thm]{Lemma}
\newtheorem{cor}[thm]{Corollary}
\newtheorem{prop}[thm]{Proposition}

\theoremstyle{definition}

\newtheorem{defin}[thm]{Definition}
\newtheorem{notat}[thm]{Notation}

\theoremstyle{remark}
\newtheorem{re}[thm]{Remark}

\numberwithin{equation}{section}

\let \al         =\alpha
\let \be         =\beta
\let \ga         =\gamma
\let \de         =\delta
\let \ep         =\varepsilon

\let \io         =\iota
\let \ka         =\kappa
\let \la         =\lambda
\let \si         =\sigma

\let \Ga         =\Gamma

\let \phi         =\varphi

\newcommand{\dist}{\mathop{\mathrm{dist}}\nolimits}
\newcommand{\Image}{\mathop{\mathrm{Im}}\nolimits}
\newcommand{\Ker}{\mathop{\mathrm{Ker}}\nolimits}

\newcommand{\ptn}{\mathbin{\widehat{\otimes}}}

\newcommand{\ptna}{\mathbin{\widehat{\otimes}_A}}

\newcommand{\rad}{\mathop{\mathrm{rad}}\nolimits}
\newcommand{\trad}{\mathop{\mathrm{t}\text{-}\mathrm{rad}}\nolimits}

\newcommand{\Ah}{\mathop{{_A}\mathrm{h}}\nolimits}

\newcommand{\N}{\mathbb{N}}
\newcommand{\CC}{\mathbb{C}}

\begin{document}
\baselineskip=17pt

\title{Topological radical of a Banach module}
\author[O. Aristov]{Oleg Yu. Aristov}
\address{36, 122 Lenin str.\\ 249034 Obninsk,
Russia}
\email{aristovoyu@inbox.ru}

\date{}

\begin{abstract}  
We introduce a concept of  topological radical of a Banach module. This 
 submodule is closed and has two descriptions: as the intersection of ranges of 
maximal contractive monomorphisms  (from the outside) and as the union of ranges 
of small morphisms (from the inside).  The topological radical is an analytic 
analogue of the radical of a module over a unital ring and has similar 
categorical properties.
  \end{abstract}



\maketitle  
     
 \section{Intoduction}    
A consideration of projective covers in \cite{Ar3} induce us to seek some 
analogue for the notion of small submodule in the Banach module context.  
Let us remind that a submodule $Y$ in a module $X$ over a ring is called 
\emph{small} (other nicks are  'superfluous' and 'coessential')  if for a
submodule $Z$ in $X$, $Y+Z=X$ implies $Z=X$. A generalization of   Dixon's  
theorem  on topologically nilpotent Banach algebras (see 
Theorem~\ref{stDinonf}) leads us to the definition of a \emph{small 
morphism}. The range of a small morphism  is a submodule in a Banach module 
and can be considered as a functional analytic analogue of small submodule. 

Our main aim is to extend the concept of Jacobson radical from Banach 
algebras to Banach modules.  As a pattern we take  the notion of radical of  
a module from Rings Theory.  But our approach offers some functional 
analytic modifications. The Jacobson radical  of a unital ring can be 
described as the \emph{intersection of all maximal left ideals }(from 
outside) or as the set  of  \emph{all  $r$ such that  $1+ar$ is invertible 
for every $a$} (from inside). This concept applies well to a unital Banach 
algebra $A$  because  \emph{every maximal left ideal is closed}  and  $1+ar$ 
is invertible for every $a\in A$ iff \emph{ $ar$ is topologically nilpotent  
(i.e. $\|(ar)^n\|^{1/n}\to 0$ for every $a\in A$)}.  

On the other hand, it is well known that the notion of radical can be extended to modules. The 
radical of a unital module $X$ over a unital ring is the \emph{intersection 
of all maximal submodules}  and  coincides with the \emph{union of all small submodules}  
(the notation is $\rad X$).  Note that for an element $r$ of a ring $A$, 
the submodule $Ar$ is small iff $1+ar$ is invertible for every $a$.  The pure algebraic 
notion of  radical of a module is useful in Banach Module Theory only  
in particular cases, for example, for finitely-generated modules~\cite{Ar3}.  
In general, neither a maximal submodule  nor a  small submodule  in a Banach 
module need not  be closed. But then again we can not restrict ourselves to 
some classes of closed submodules because  submodules of the form $A\cdot x$ 
(that are potentially not closed) play an important role in the basic theory 
of module radicals.  As we see below the right way is to consider ranges of 
bounded module morphisms as an intermediate class between closed submodules 
and all submodules. But it is seems more appropriate from the ideological 
and the technical points of view to work with  morphisms themselves instead of 
their ranges.  

In this article we introduce a concept of  \emph{topological radical} of a 
Banach module. This closed submodule  has two descriptions: as the 
intersection of ranges of maximal contractive monomorphisms  (from the outside) 
or as the union of ranges of small morphisms (from the inside). 

The author would like to thank the referee for the valuable comments which 
helped to improve the manuscript.
     
\section{Small morphisms of Banach modules}   

Let $A$ be a Banach algebra. We suppose that the norm of the multiplication in 
$A$ is not greater than $1$. For $n\in\N$ set
\begin{equation}\label{topnild}
S(n)\!:=\sup\|r_1r_2\cdots r_n\|^{1/n},
\end{equation}
where $r_1,r_2,\ldots, r_n$ run over the unit ball of $A$.  If $\lim_{n\to 
\infty} S(n)=0$ then $A$ is called \emph{topologically nilpotent}. Note that 
$A$ is topologically nilpotent if and only if for every bounded sequence 
$(r_n)\subset A$, $$ \lim_{n\to\infty}\|r_1r_2\cdots r_n\|^{1/n}=0. $$ 
Obviously, a topologically nilpotent Banach algebra is radical.

Recall that $C[0,1]$ and $L^1[0,1]$ are radical Banach algebras with respect to the 
cut-off convolution. The first algebra  is topologically nilpotent but the 
second algebra is not topologically nilpotent \cite[Section 4.8.8]{Pa}.

In \cite{Di} P.~G. Dixon shows that $A\cdot X\neq X$ for every non-trivial 
left Banach module $X$ over a topologically nilpotent Banach algebra $A$  
(see the proof in \cite[Theorem 4.8.9]{Pa} also). But in fact his argument 
gives a stronger assertion. Let us set $$ \pi^A_X\!:A\ptn X\to X\!:r\otimes 
x\mapsto r\cdot x $$  for a left Banach $A$-module $X$, where $\ptn$ denotes 
the projective tensor product of Banach spaces. (We suppose that the norm of 
multiplication in $X$ is not greater than $1$ also.)  Below "a module" means an $A$-module.
\begin{thm}[Dixon]\label{Dinonf}
If $X$ is a non-trivial left Banach module  over a topologically nilpotent 
Banach algebra $A$, then $\Image \pi^A_X\neq X$.
\end{thm}
We need a more general result.
\begin{thm}\label{stDinonf}
Let $A$ be a topologically nilpotent Banach algebra, and let $\phi\!:Y\to X$ 
be a morphism of left Banach modules such that $X=\Image \phi+\Image 
\pi^A_X$. Then $\phi$ is surjective.
\end{thm}
\begin{proof}
The assumption of the theorem means that the morphism $$Y\oplus (A\ptn X) 
\to X\!: (y,u)\mapsto \phi(y)+\pi^A_X(u) $$ is surjective. (Here the sum is 
endowed with the $\ell^1$-norm.)  By the open mapping theorem there is $C>0$ 
with the following property. For every $x$ in $X$ there exist $y\in Y$, 
$r_i\in A$, and $x_i\in X$ such that 
\begin{equation}
\label{openmap}
x=\phi(y)+\sum_{i=1}^\infty r_i\cdot 
x_i\qquad\text{and}\qquad \|y\|+\sum_i \|r_i\|\,\|x_i\|\le C\|x\|.
\end{equation} 

Now we fix $x$ in $X$ and choose by induction sequences $(y_n)\subset Y$ and 
$(v_n)\subset X$  such that
\begin{equation}
\label{xxx0} x=\phi(y_n)+v_n, 
\end{equation} 
where $v_n$ can be represented as
\begin{equation}
\label{xxx1}
 v_n=\sum_{i=1}^\infty r_{1,i}\cdots r_{n,i} \cdot x_i\quad
 \text{for some $r_{1,i},\ldots,r_{n,i}\in A$ and $x_i\in X$; $i\in\N$,}
\end{equation}
and the following two conditions are satisfied.
\begin{equation}\label{xxx3}
 \|y_{n+1}-y_n\|\le C \sum_i \|r_{1,i}\cdots r_{n,i}\|\,\|x_i\|;
\end{equation}
\begin{equation}\label{xxx2}
\sum_i \|r_{1,i}\|\cdots \|r_{n,i}\|\,\|x_i\|\le
 C^n\|x\|.
\end{equation}

Suppose that for $n\in\N$ we have elements $y_1,\ldots,y_n$ and 
$v_1,\ldots,v_n$ that satisfy the above conditions, in particular, the 
condition~(\ref{xxx3}) satisfies up to $n-1$. Fix decompositions in 
(\ref{xxx0}) and (\ref{xxx1}). Set $t_i\!:=r_{1,i}\cdots r_{n,i}$. Applying 
(\ref{openmap}) we can write every $x_i$ as
\begin{equation}\label{xxx4}
x_i=\phi(y_i')+\sum_j s_{ji}\cdot 
x_{ji}',\quad\text{where}\quad\|y_i'\|+\sum_j \|s_{ji}\|\,\|x_{ji}'\|\le 
C\|x_i\|.
\end{equation}
Then $$ x=\phi(y_n)+v_n=\phi(y_n)+\sum_{i}t_i\cdot\phi(y_i')+\sum_{i,j} t_i 
s_{ji}\cdot x_{ji}'. $$ Now set $ y_{n+1}\!:=y_n+\sum_{i} t_i\cdot y_i'$ and 
$v_{n+1}\!:=\sum_{i,j} t_i s_{ji}\cdot x_{ji}'$. It follows 
from~(\ref{xxx4}) that $\|y_i'\|\le C\|x_i\|$. Hence, $$ \|y_{n+1}-y_n\|\le 
\sum_{i} \|t_i\|\,\|y_i'\|\le C \sum_{i} \|t_i\|\,\|x_i\|, $$ i.e. we 
obtain~(\ref{xxx3}). By~(\ref{xxx2}) and~(\ref{xxx4}) we get 
$$
 \sum_{i,j} 
\|r_{1,i}\|\cdots \|r_{n,i}\|\,\|s_{ji}\|\,\|x'_{ji}\|\le C^{n+1}\|x\|, 
$$ i.e. 
after an obvious change of notation we have~(\ref{xxx1}) and~(\ref{xxx2}) 
for $n+1$. By induction, there exist sequences with the desired properties.

Note that (\ref{xxx2}) implies (see (\ref{topnild})) $$ \sum_i \|r_{1,i}\cdots r_{n,i}\|\,\|x_i\|\le 
\sum_i S(n)^{n}\|r_{1,i}\|\cdots \|r_{n,i}\|\,\|x_i\|\le S(n)^nC^n\|x\| $$ for 
every $n$.  Therefore $\|v_n\|\le S(n)^nC^n\|x\|$ and  $\|y_{n+1}-y_n\|\le 
S(n)^nC^{n+1}\|x\|$ by (\ref{xxx3}). Hence we have for $m>n$
\begin{equation*}
\|y_{m}-y_n\|\le  \sum_{k=n}^{m-1} S(k)^kC^{k+1}\|x\|.
\end{equation*}
Since $ S(n)\to 0$, it follows that $y_n$ is a fundamental sequence and 
$v_n\to 0$. Finally, from $x=\phi(y_n)+v_n$ we get $x=\phi(\lim_n y_n)$, 
i.e. $x\in\Image\phi$.
\end{proof}       
 
\begin{notat}\label{dotpl}  
Let $\phi\!:Y\to X$ and $\psi\!:Z\to X$ be morphisms of Banach modules. 
Denote by $\phi\dotplus\psi$ the morphism
$$
Y\oplus Z\to X\!:(y,z)\mapsto \phi(y)+\psi(z).
$$ 

\end{notat}

\begin{defin}\label{smallmor}
We say that a morphism $\psi\!:X_0\to X$ of Banach modules is \emph{small} 
if for every morphism $\phi\!:Y\to X$ such that $\phi\dotplus\psi$ is 
surjective  $\phi$ is surjective  also, i.e., $\Image \phi+ \Image \psi = X$ implies $\Image \phi = X$.
\end{defin}

Thus, Theorem~\ref{stDinonf} asserts that \emph{for every  left Banach 
module $X$ over a topologically nilpotent Banach algebra $A$ the morphism 
$\pi^A_X$ is small}.    

\begin{prop}\label{compsm} 
If $\psi\!:X_0\to X$ is a small morphism then $\tau\psi$ is small for  each 
module $V$ and each morphism $\tau\!:X \to V$.
\end{prop}     
\begin{proof}
Let $\tau\!:X\to V$ and $\phi\!:Y\to V$  be morphisms of Banach modules such 
that  $\phi\dotplus\tau\psi$ is surjective. Consider the pullback diagram  
 \begin{equation*}% \label{}
 \xymatrix@C=20pt{  Y\times_V X
\ar[r]\ar[d]_{\phi'} &  Y  \ar[d]^\phi \\
   X   \ar[r]_\tau &  V  }    
 \end {equation*}                                                        
associated with $\phi$ and $\tau$.   

For every $x\in X$ there are $y\in Y$ and $z\in X_0$ such that 
$\tau(x)=\phi(y)+\tau\psi(z)$. Then $\phi(y)=\tau(x-\psi(z))$. By explicit 
construction of $Y\times_V X$ this means that   $w=(y,x-\psi(z))\in 
Y\times_X V $ and $\phi'(w)= x-\psi(z)$. Hence, $ \phi'\dotplus\psi$ is 
surjective.  Since $\psi$ is small,    $\phi'$ is  surjective also.   
Therefore for every $x\in X$ there exists $y\in Y$ such that 
$\phi(y)=\tau(x)$.    The assumption that $\phi\dotplus\tau\psi$ is 
surjective implies that $\phi\dotplus\tau$ is surjective. Thus, $\phi$ is  
surjective also.    
\end{proof}

\begin{prop}\label{epsm}
Let $\psi\!:X_1\to X$ be a morphism  of Banach modules, and let 
$\ep\!:X_0\to X_1$ be a surjective morphism  of Banach modules such that 
$\psi\ep$ is small. Then $\psi$ is small.
\end{prop}
\begin{proof}
Suppose that $\phi\!:Y\to X$ is a morphism such that $\phi\dotplus\psi$ is 
surjective.
 Then $\phi\dotplus\psi\ep$ is surjective
also. Since $\psi\ep$ is small, $\phi$ is surjective.
\end{proof}

Recall that a left Banach $A$-module $P$ is called  \emph{strictly 
projective} if for each  surjective morphism of Banach $A$-modules 
$\ep\!:Y\to P$ there exists a morphism $\rho\!:P\to Y$ such that 
$\ep\rho=1$. Denote by $\ell^1$ the infinite Banach $\ell^1$-space with a 
countable basis.

\begin{thm}\label{aaitn}
{ \rm (cf.~\cite[Th.11.5.5]{Ka})} Let $I$ be a closed left ideal in a unital 
Banach algebra $A$, and let $\io\!:I\to A$ be the natural inclusion. The 
following conditions are equivalent.

{\em (A)}  $I$ is topologically nilpotent.

{\em (B)} For every unital left Banach $A$-module $X$ the morphism of
Banach $A$-modules $I\ptna X\to X\!:a\otimes_A x\mapsto a\cdot x$ is small.

{\em (C)} For every strictly projective unital left Banach $A$-module $P$ the morphism of
Banach $A$-modules $I\ptna P\to P\!:a\otimes_A x\mapsto a\cdot x$ is small. 

{\em (D)} The morphism of left Banach $A$-modules $ (\io\otimes
1)\!:I\ptn\ell^1\to A\ptn\ell^1 $ is small.
\end{thm}    
 \begin{proof}
(A)$\Rightarrow$(B)
 If $I$ is topologically nilpotent and $X$ is a left Banach $A$-module then by
Theorem~\ref{stDinonf} $\pi^I_X$ is small as a morphism of left Banach 
$I$-modules. Hence, it is small as a morphism of left Banach $A$-modules. 
Since $\pi^I_X$  is a composition of a surjective morphism $I\ptn X\to 
I\ptna X$ and a morphism $I\ptna X\to X$, Proposition~\ref{epsm} implies~(B).
                                                                               
(B)$\Rightarrow$(C) It is obvious.

(C)$\Rightarrow$(D) It is easy to see that $A\ptn\ell^1$ is strictly 
projective. By assumption
\begin{equation}\label{ail1}
I\ptna A\ptn \ell^1\to A\ptn \ell^1 \!:a\otimes_A b \otimes x\mapsto 
ab\otimes x
\end{equation}
is a small morphism of left Banach $A$-modules. Since $A$ is unital, $I\ptna 
A\cong I$, and we have~(D).

(D)$\Rightarrow$(A). Let $(a_n)$ be a bounded sequence in $I$, and let 
$\{e_i\}_{i\in\N}$ be the canonical basis in $\ell^1$. Consider
\begin{equation}\label{submy}
\phi\!:A\ptn\ell^1\to A\ptn\ell^1\!:\sum_{i=1}^\infty b_i\otimes e_i\mapsto 
\sum_{i=1}^\infty b_ia_i\otimes e_{i+1}.
\end{equation}
It is obvious that $\phi$ is a morphism of left Banach modules. Fix $\la\in 
\CC$. Since $$ \sum_{i=1}^\infty  b_i a_i\otimes e_{i+1}\in I\ptn \ell^1, $$ 
we have $I\ptn \ell^1+\Image(1+\la\phi)=A\ptn\ell^1$. Since$(\io\otimes 1)$ 
is small, $1+\la\phi$ is surjective. If $(1+\la\phi)(u)=0$ for some 
$u=\sum_i b_i\otimes e_i$, then $b_1=0$ and $b_{i+1}-\la b_ia_i=0$ for all 
$i$. It follows that $b_i=0$ for all $i$, so that $1+\la\phi$ is injective. 
Thus, $1+\la\phi$ is an isomorphism for every $\la\in\CC$. This implies that 
$\phi$ is a topologically nilpotent operator, i.e. 
$\lim_{n\to\infty}\|\phi^n\|^{1/n}=0$.

It is clear that $$ \|\phi^n(1\otimes e_1)\|=\|a_1a_2\cdots a_n \otimes 
e_{n+1}\|=\|a_1a_2\cdots a_n\|. $$ Therefore $ \|a_1a_2\cdots a_n\|\le 
\|\phi^n\|$. The rest is obvious.
\end{proof} 

Considering every Banach algebra as an ideal in the unitization we have

\begin{cor}\label{chartnba}
A Banach algebra $A$ is topologically nilpotent if and only if for every 
Banach $A$-module $X$ the morphism $A\ptn_A X\to X$ is small if and only if  
for every strictly projective  left Banach $A$-module $P$ the morphism 
$A\ptn_A P\to P$ is small. 
\end{cor}

Note that the definition of $S(n)$ is invariant under replacement of the 
left multiplication by the right multiplication. So all results above can be 
applied to right Banach modules. 

If $X$ and $Y$ are left Banach $A$-modules  we denote by $\Ah(X,Y)$  the set 
of all bounded  $A$-module morphisms  from $X$ to $Y$. Recall that a left 
$A$-module $X$ is called \emph{unital}  if $1\cdot x = x$ for all $x\in X$.

\begin{prop}\label{alsmXY}
Let $A$ be a Banach algebra, $X$ and $Y$  unital left Banach $A$-modules, 
and  $\al$ in $\Ah(X,Y)$. The following conditions are equivalent.  

{\em (1)}    $1-\al\phi$ is right invertible in the unital algebra 
$\Ah(Y)$   for every $\phi \in \Ah(Y,X)$.

{\em (2)}    $\al \circ\Ah(Y,X)$ is a small right ideal in  $\Ah(Y)$.

\end{prop}
\begin{proof}
$(1)\Rightarrow(2)$ Let $L$ be a  right ideal in $\Ah(Y)$  such that $\al 
\circ\Ah(Y,X)+L=\Ah(Y)$. Then there are $\phi \in \Ah(Y,X)$  and $\psi\in L$ 
satisfying $\al\phi+\psi=1$. By assumption $\psi$ has a right inverse 
$\psi_1$, hence, as $L$ is a right ideal  in $\Ah(Y)$, we have $1 = 
\psi\psi_1\in L$, so that $L = \Ah(Y )$.

$(2)\Rightarrow (1)$ Let $\phi \in \Ah(Y,X)$. Set 
$L\!:=(1-\al\phi)\circ\Ah(Y)$. Then $1-\al\phi\in L$; so that  $1\in \al 
\circ\Ah(Y,X)+L$. Therefore $\al \circ\Ah(Y,X)+L=\Ah(Y)$. Since $\al 
\circ\Ah(Y,X)$ is small, $L=\Ah(Y)$. This implies that $1-\al\phi$ is right 
invertible.
\end{proof}
 
\begin{thm}
Let $X$ and $P$  be  unital left Banach $A$-modules. Suppose that $P$ is  
strictly projective and $\al\in\Ah(X,P)$. The following conditions are 
equivalent.

{\em (1)} $\al$ is small.   

{\em (2)} $1-\al\phi$ is right invertible in $\Ah(P)$   for all $\phi \in 
\Ah(P,X)$ .  

{\em (3)} $\al\circ\Ah(P,X)$ is a small right $\Ah(P)$-submodule in $\Ah(P)$.
\end{thm}
\begin{proof}
$(1)\Rightarrow (2)$ Let $\phi \in \Ah(P,X)$. Then $(1-\al\phi)\dotplus 
\al\!:P\oplus X\to P$ is obviously surjective. Since $\al$  is small, 
$1-\al\phi$ is surjective.  Since $P$ is strictly projective, $1-\al\phi$ 
admits a right inverse morphism.

$(2)\Rightarrow(1)$ Suppose $\eta\in\Ah(Y,P)$  for some $Y$ and 
$\eta\dotplus\al$ is surjective. Since $P$ is strictly projective, 
$\eta\dotplus\al $ is right invertible, i.e.  there exist $\psi_1\in 
\Ah(P,Y)$ and $\psi_2\in \Ah(P,X)$  such that $\eta\psi_1+\al\psi_2=1$. By 
assumption $\eta\psi_1$ is right invertible. Hence, $\eta$ is surjective.

$(3)\Leftrightarrow(2)$ follows from Proposition~ \ref{alsmXY}.
\end{proof}     

A surjective morphism $\ep\!:X\to V$ of Banach $A$-modules is said to be a 
\emph{cover} if a morphism $\phi\!:Y\to X$ of Banach $A$-modules is a 
surjective morphism, whenever $\ep\phi$ is so \cite{Ar3}.

\begin{prop}\label{covsmall}
A surjective morphism  $\ep\!:X\to V $ of Banach modules is a cover  if and 
only if the embedding $\Ker \ep\to X$ is a small morphism.
\end{prop}
\begin{proof}  
Denote the embedding $\Ker \ep\to X$ by $\ker\ep$. Suppose that $\ep$ is a 
cover. Let $\phi\!:Y\to X$ be a morphism of Banach modules such that 
$\phi\dotplus\ker\ep$ is surjective. Note that $\ep(\phi\dotplus\ker\ep)= 
\ep\phi $  is surjective also.   Since $\ep$ is a cover, $\phi$ is 
surjective. Thus, $\ker\ep$ is a small morphism.

Suppose that $\ker\ep$ is small. Let $\phi\!:Y\to X$ be a morphism of Banach 
modules such that $\ep\phi$ is surjective. Then for every $x\in X$ there is 
$y\in Y$ such that $\ep(x)=\ep\phi(y)$. Then  $x=\phi(y)+(x-\phi(y))$ where 
$ x-\phi(y)\in \Ker\ep$. Therefore $\phi\dotplus\ker\ep$ is surjective. 
Since $\ker\ep$ is small, $\phi$ is surjective. Thus, $\ep$ is a cover.
\end{proof}

 \section{Maximal contractive monomorphisms}    
 \label{smax}
 
 Fix a unital Banach algebra $A$ and a left unital Banach $A$-module $X$. 
Consider a pre-order on the  set of contractive monomorphisms with range in  
$X$. We set $\be\succeq\ga$  for  $\be$ and $\ga$ if  there exists a 
contractive morphism $\ka$ such that $\ga=\be\ka$. We say that $\be$ and 
$\ga$ are equivalent if  $\ka$ is an isometric isomorphism. The pre-order 
induces an order on the set of equivalence classes of contractive 
monomorphisms. 

\begin{re}
If $X$ is unital and $\be\!:Y\to X$ is a monomorphism then $Y$ is also 
unital. To see this consider the decomposition $Y=Y_0\oplus Y_1$,  where 
$Y_0=\{y\in Y:\,1\cdot y=0 \}$  and $Y_1=\{y\in Y:\,1\cdot y=y  \}$. Since 
$X$ is unital, $\Ah(Y_0,X)=0$. Therefore $Y_0=0$.  Thus, we does not need  
the restriction on initial module of a monomorphism.

\end{re} 

\begin{defin} \label{meetjoinm} 
     Let $\be\!:Y\to X$ and $\ga\!:Z\to X$  be contractive monomorphisms.
     
     (1)  Denote by $\be\vee \ga$ the natural morphism  $(Y\oplus Z)/\Ker(\be\dotplus\ga)\to X$ 
     associated with $\be\dotplus\ga$.
     
     (2)  Denote by $\be\wedge \ga$ the natural morphism $Y\times_X Z \to X$, where $Y\times_X Z$ is
     the pullback of $\be$ and $\ga$.
\end{defin}

It is not hard to check that  $\be\vee \ga$  and $\be\wedge \ga$ are   
contractive monomorphisms. For equivalence classes $[\be]$ and $[\ga]$  of 
$\be$ and $\ga$  we set 
 $[\be]\vee [\ga]\!:=[\be\vee \ga]$  and   $[\be]\wedge [\ga]\!:=[\be\wedge \ga]$.  It is easy to see that
 these operations are well-defined.      
 
 \begin{prop}\label{supinfm}  
 Let $\be$ and $\ga$  be contractive monomorphisms.
 Then, with respect to the order define above, $[\be]\vee [\ga]$   and   $[\be]\wedge [\ga]$  
  are the supremum and the infimum  of $[\be]$ and $[\ga]$,
 respectively.        
 \end{prop}
The proof is standard. 

\begin{defin} \label{maxmono}
We say that a contractive monomorphism of left unital Banach $A$-modules 
$\al\!:Y\to X$ is \emph{maximal} if  $\al$ is not surjective and for every 
non-surjective contractive monomorphism $\be$ and every contractive morphism 
$\ka$ the equality $\al=\be\ka$ implies that $\ka$ is an isometric 
isomorphism. 
\end{defin}

 Thus, $\al$ is maximal iff
$[\al]$  is maximal  in the set  of equivalence  classes of all 
non-surjective monomorphisms  with range in $X$.   

Recall that a morphism $\ep\!:Y\to X$ is called a $C$-epimorphism for some 
$C\ge1$  if for every $x\in X$ there exist $y\in Y$  such that
 $x=\phi(y)$ and
 $ \|y\|\le C\|x\|$.                                        
  
\begin{prop}\label{1oc} 
For $x_0\in X$, set
$$
  \tau\!:A\to X\!:a\to a\cdot x_0.
$$    
Suppose that  $\phi\!:Y\to X$ is a morphism  such that $x_0\notin\Image\phi$ 
and
 $\phi\dotplus\tau$  is a $C$-epimorphism for $C\ge1$. Then $\dist(x_0,\Image\phi)\ge 1/C$.                
 \end{prop}
\begin{proof} 
Assume that  
 $\|x_0-\phi(y)\|<1/C$ for some $y\in Y$.  
Since $\phi\dotplus\tau$ is $C$-epimorphism, there exist $y'\in Y$ and $a\in A$ such that
 $x_0-\phi(y)=\phi(y')+a\cdot x_0$ and
 $$
 \|y'\|+\|a\|\le C\|x_0-\phi(y)\|<1.
 $$       
  
 Thus, $\|a\|<1$, hence  $1-a$ is invertible in $A$.     Therefore
   $$
 x_0=\phi((1-a)^{-1}\cdot(y'+y)).
 $$ 
Hence,  $x_0\in \Image \phi$.  We get a contradiction.  
\end{proof}  

\begin{thm}\label{maximcl}
  Every maximal contractive monomorphism is an isometry.
\end{thm}  
\begin{proof}    
Let $\al\!:Y\to X$ be a maximal contractive monomorphism, so $\Image \al\ne 
X$. Suppose that $x_0\in X\setminus\Image \al$.  Define $\tau$ as in 
Proposition~\ref{1oc}. Denote by $\be$ the natural monomorphism $(Y\oplus 
A)/\Ker(\al\dotplus\tau)\to X$ and by $\ka$ the composition $$Y\to Y\oplus 
A\to  (Y\oplus A)/\Ker(\al\dotplus\tau) .$$

We claim that $\be$ is surjective. Indeed, assume the converse.  Since $\al$ 
is maximal, $\ka$ is an isometric isomorphism. In particular, there is $y\in 
Y$ such that
$$(y,0)-(0,1)\in \Ker(\al\dotplus\tau).$$ Hence, $\al(y)=x_0$. We get a contradiction.

Since $\be$ is surjective,  $ \al\dotplus\tau$ is surjective.  By the open 
mapping theorem,  $\al\dotplus\tau$ is a $C$-epimorphism for some $C\ge 1$.  
It follows Proposition~\ref{1oc} that $\dist(x_0,\Image\phi)\ge 1/C$. Since 
$x_0$ is arbirtary,  $\Image \al$ is closed. Let $\ga\!:\Image\al\to X$  be 
the natural embedding. Since $\al=\ga\al$ and $\al$ is maximal, $\al$ is an 
isometry. 
\end{proof} 

\begin{thm}\label{maxirr}
  Let $X'$ be a closed submodule of $X$. Then the natural embedding $\io\!:X'\to X$ 
  is a maximal contractive monomorphism if and only if $X/X'$ is an irreducible module.
\end{thm}  
\begin{proof}    
 $(\Rightarrow)$  Assume that $\io$ is maximal. Let $x_0\in X\setminus X'$ and $x_1\in X$.    Since $\al$ is maximal,
 $X'+A\cdot x_0=X$. In particular, there is $a\in A$ such that $x_1-a\cdot x_0\in X'$. Therefore $x_0+X'$ is
 a cyclic element of $X/X'$. Hence,  $X/X'$ is irreducible.
 
 $(\Leftarrow)$ Assume that $X/X'$ is irreducible. Suppose that there are a non-surjective contractive monomorphism $\be$ and
 a contractive morphism $\ka$ such that $\be\ka=\io$.  
 
 Since $X'\subset \Image\be$, $\Image\be\neq X$ and $X/X'$ is irreducible, 
  $\Image\be= X'$. Therefore, $\be\ka=1$. Since $\be$ is a monomorphism, it is an isomorphism.  
  Since $\be$ and $\ka$ are contractive, $\ka$ is isometric.
\end{proof}        

Note that $X/X'$ is irreducible iff $X'$ is a maximal submodule in the 
algebraic sense. Thus, maximal monomorphisms can be described as embeddings 
of closed  maximal submodules.     
       
% Below we need the lifting Lemma for maximal monomorphisms.        

\begin{lem}\label{lifmax}
Let $Z$ be a closed submodule in $X$ and  $\al\!:Y\to X/Z$  a maximal 
contractive monomorphism. Denote  the projection $X\to X/Z$ by $\si$. Then 
there exists a commutative diagram 

\begin{equation*}% \label{}
 \xymatrix@C=20pt{ W
\ar[r]^\mu\ar[d]_\be &  Y  \ar[d]^{\al} \\
   X   \ar[r]^{\si} &  X/Z  }
\end{equation*}  
where $\be$ is a maximal contractive monomorphism.
\end{lem}
\begin{proof} 
Set $W\!:=Y\times_{X/Z}X$ or, more precisely, $W=\{(y,x)\in Y\times X\!: 
\al(y)=\si(x)\}$. Denote by $\be$ and $\mu$ the morphisms $(y,x)\mapsto x$ 
and  $ (y,x)\mapsto y$, respectively. Note that $\mu$ is surjective and 
$Z\subset \Image \be$.    It is obvious that $\be$ is a contractive 
monomorphisms.

Suppose that $\ga\!:V\to X$ is a non-surjective contractive monomorphism  
and  $\ka\!:W \to V$ is a contractive morphism such that $\be=\ga\ka$.  

Assume that  $\Image\al+\Image(\si\ga)=X/Z$. Since $\Image \al=\Image 
(\al\mu)=\Image(\si\be)=\Image(\si\ga\ka)$, we have $\Image(\si\ga)=X/Z$.    
It follows from $Z\subset\Image \be \subset\Image \ga$ that $\ga$ is 
surjective. We get a contradiction. Hence, $\Image\al+\Image(\si\ga)\ne 
X/Z$. 

Since $\al$ is maximal, $\Image(\si\ga)\subset\Image\al$. By 
Theorem~\ref{maximcl} $\al$ is an isometry,  therefore there is a 
well-defined contractive morphism $\de:\!V\to Y$ such that $\al\de=\si\ga$. 
The pull-back property implies that there is a contractive morphism 
$\rho\!:V\to W$ such that  $\be\rho=\ga$. Then  $[\be]=[\ga]$.  Thus, $\be$ 
is maximal.
\end{proof}     
   
 \begin{lem}\label{smmax0}     
Let  $\al$ and $\be$ be non-surjective contractive monomorphism with ranges 
in $X$ such that 
  $\be\succeq\al$, and let $\phi$ be a morphism  such that $\al\dotplus\phi$ is a $C$-epimorphism for some $C\ge1$.
  Then $\be\dotplus\phi$ is a $C$-epimorphism.
   \end{lem}
   \begin{proof}
   Suppose that   $\ka\!:Y\to Z$ is a contractive morphism such that
$\al=\be\ka$.      Since $\al\dotplus\phi$ is a $C$-epimorphism, for every 
$x\in X$ there exist $x_0\in X_0$ and $y\in Y$ such that 
$x=\phi(x_0)+\al(y)$ and $\|x_0\|+\|y\|\le C\|x\|$. Denote $\ka(y)$ by $z$. 
Then    $x=\phi(x_0)+\be(z)$ and $\|x_0\|+\|z\|\le \|x_0\|+\|y\|\le C\|x\|$. 
Therefore, $\phi\dotplus\be$ is a $C$-epimorphism.
   \end{proof} 
 
 \begin{lem}\label{smmax1} 
 Let  $C\ge1$, and let $\phi$ be a contractive
morphism with range in $X$. Denote by $\Ga$ a family of all  contractive 
monomorphisms $\al$  with range in $X$
 such that
  
{\em (1)}   $\al$  is not surjective;

{\em (2)} $\al\dotplus\phi$ is a $C$-epimorphism.
 
Suppose that there are $\de>0$ and $x_0\in X$ such that  $\dist(x_0,\Image\al)\ge\de$
for every $\al\in \Ga$.
Then for every $\al_0\in\Ga$  there exists a maximal contractive
monomorphism $\ga$ such that $\ga\in\Ga$ and
 $\ga\succeq\al_0$. 
 \end{lem}  

 \begin{proof} 
Set $\Ga'\!:=\{\al\in \Ga; \al\ge\al_0\}$.  
Suppose that $\Ga_0$ is a linear ordered subset of $\Ga'$. We claim that that $\Ga_0$ admits an upper bound.
 
 Denote by $Y_\al$ the initial module of $\al\in\Ga_0$ and by $\ka_{\al\al'}$ the connecting contractive morphism
 for $\al$ and $\al'$ in $\Ga_0$ such that $\al'\succeq\al$.
 Then there exists an inductive limit $Y$ of a spectral family $(\ka_{\al\al'})$ 
 in the category of contractive  morphisms. In particular,  
 there is a family $(\ka_{\al}\!:Y_\al\to Y)$ of contractive  morphisms and 
  $\be\!:Y\to X$ such that $\al=\be\ka_{\al}$ for every $\al$. 
  Note that $\bigcup_{\Ga_0}\Image\ka_{\al}$ is dense in $Y$, hence
  $\bigcup_{\Ga_0}\Image\al$ is dense in $\Image\be$. 
   Since
  $\dist(x_0,\Image\al) \ge\de$ for all $\al$, we have $\dist(x_0,\Image\be)\ge \de$.  Hence,
  $\be$ is    not surjective.
     Applying  Lemma~\ref{smmax0} 
   we get that            
  $\be\dotplus\phi$ is a $C$-epimorphism.  Therefore,
  we have that $\be\in\Ga'$ and $\be\succeq\al$ 
 for every $\al \in \Ga_0$. 
 
 Since $\Ga'$ is not empty and every linear ordered subset  in 
 $\Ga'$ admits an upper bound, 
 there is a maximal element $\ga$ in $\Ga'$. 
 Now we claim that $\ga$ is a maximal contractive monomorphism.      
 Suppose that $\ga'\!:Z\to X$ is a non-surjective contractive monomorphism  
and  $\ka\!:Y\to Z$ is a contractive morphism such that
$\ga=\ga'\ka$.  It follows from Lemma~\ref{smmax0} that 
$\ga'\dotplus\phi$ is a $C$-epimorphism. Hence, $\ga'\in \Ga'$. Since $\ga$  
is maximal in $\Ga'$, $\ka$ is an isometric isomorphism. 
Thus  $\ga$ is a maximal contractive monomorphism. By construction $\ga\succeq\al_0$. 
  \end{proof}  
       
\begin{re}
In the proof we find the supremum of a directed set of contractive 
monomorphism implicitly.   It is not hard to see that the constructions  of  
$\vee$ and $\wedge$  from Definition~\ref{meetjoinm}  can be applied to 
arbitrary sets of  monomorphism also.
\end{re}

\begin{prop}\label{fdmax} 
Suppose that $X$ is finitely-generated. Then for every non-surjective 
contractive monomorphism $\al_0$  with range in $X$ there exists a maximal 
contractive monomorphism $\ga$ such that $\ga\succeq\al_0$.
 \end{prop}
\begin{proof}
Let $x_1,\ldots,x_n$ be generators of $X$. Consider morphisms
   $$\tau_i\!:A\to X\!:a\to a\cdot x_i\quad (i=1,\ldots,n).$$   
  Since $\al_0$ is not surjective and $\al_0\dotplus\tau_1\dotplus\cdots\dotplus\tau_n$ is surjective,
  there exist minimal $k$ in $\{1,\ldots,n\}$ 
such that  $\al_0\dotplus\tau_1\dotplus\cdots\dotplus\tau_k$ is surjective. 
This implies that
        there exists $C\ge 1$ such that
    $\al_0\dotplus\tau_1\dotplus\cdots\dotplus\tau_k$   is a $C$-epimorphism. 
    
    Denote by $\Ga$ a family of all  contractive
monomorphisms $\be$  with range in $X$
 such that $\be$  is not surjective and
$\be\dotplus\tau_k$ is a $C$-epimorphism. Let us set 
$\be_0=\al_0\vee\al_1\vee\cdots\vee\al_{k-1}$, where    $\al_i$ is a 
contractive monomorphism such that
 $\Image\al_i=\Image\tau_i$. 
Note that $x_k\notin\Image\be$ for every $\be\in \Ga$.     
It follows from Proposition~\ref{1oc} that   $\dist(x_k,\Image\be)\ge 1/C$ for every $\be\in \Ga$.
Thus, the conditions of Lemma~\ref{smmax1} are satisfied. Hence,  there exists a maximal contractive
monomorphism $\ga$ such that 
 $\ga\succeq\be_0$.  Therefore  $\ga\succeq\al_0$.      
\end{proof}

\section{Topological radical of a Banach module}

Note that equivalence  classes of contractive morphisms form a lattice with 
respect to the operations $\vee$  and $\wedge$. Under some conditions there 
is a standard way to define a radical in a lattice using small and maximal 
elements (see, for example, \cite[Ch.9, Exersises]{Ka}). But on this way we 
meet two difficulties. First, we define small and maximal morphism in 
different categories  of Banach modules (the topological and the metric 
categories). Second, there are no sufficiently many compact elements in our 
lattice. However, using basic Proposition~\ref{1oc}  and its corollaries 
from Section~\ref{smax} we can find the desired topological interplay 
between small and maximal morphisms.

\begin{prop}\label{smmax}
Let $X$ be a left unital Banach $A$-module and let $x_0\in X$. If 
$$
  \tau\!:A\to X\!:a\to a\cdot x_0
$$  
is not small then there exists a maximal contractive monomorphism $\ga$ 
such that $x_0\notin\Image \ga$.
\end{prop}  
\begin{proof} 
Since $\tau$ is not small,  there exists a non-surjective  morphism 
$\al_0\!:Y\to X$ such that $\al_0\dotplus\tau$ is surjective. By  the open 
mapping theorem,  there is $C\ge 1$ 
 such that $\al_0\dotplus\tau$ is a $C$-epimorphism.  
 We can assume  $\al_0$  is a contractive monomorphism.      
 
 Denote by $\Ga$ the family of all  contractive
monomorphisms $\al$  with range in $X$
 such that   $\al$  is not surjective and
 $\al\dotplus\tau$ is a $C$-epimorphism for $C$ chosen above. Note that $x_0\notin\Image \al$ for every $\al\in\Ga$.
Proposition~\ref{1oc} implies that   $\dist(x_0,\Image\al)\ge 1/C$ for every $\al\in \Ga$.
It  follows from Lemma~\ref{smmax1} that there exists a maximal contractive
monomorphism $\ga$ such that $\ga\in\Ga$.  Hence, $x_0\notin\Image \ga$.
\end{proof}   
 
\begin{thm}\label{exsrad}
Let $X$ be a left unital Banach $A$-module. Denote  by $X_1$ the set 
$\bigcup\Image\psi$, where $\psi$ runs  all small morphism with the range in 
$X$,  and denote by $X_2$ the set $ \bigcap \Image\ga$, where $\ga$ runs all 
maximal contractive monomorphisms with the range in $X$. Then $X_1=X_2$ and 
this submodule of $X$ is closed.
\end{thm}     
\begin{proof}

{\rm (1)} Suppose that $x_0\in X_2$.   Assume that $  \tau\!:A\to X\!:a\to 
a\cdot x_0$  is not small. By Proposition~\ref{smmax}  there exists a 
maximal contractive monomorphism $\ga$ such that $x_0\notin\Image \ga$. 
Therefore $x_0\notin X_2 $. This contradiction implies that $\tau$ is small. 
Thus, $X_2\subset X_1$.

{\rm (2)}    Suppose that $\psi$ is a small morphism with the range in $X$. 
We can assume  $\psi$  is a contractive monomorphism.
Suppose that there is a maximal contractive monomorphisms $\ga$ such that $\Image\psi$  is not
a subset of  $\Image\ga$.  Then $\ga\vee\psi=1$. Therefore,
  $\Image\ga+\Image\psi=\Image(\ga\vee\psi)=X$. Since $\psi$ is small, $\ga$ is surjective. 
This contradiction implies that $\Image\psi\subset\Image\ga$.  Thus, $X_1\subset X_2$.

It follows from Theorem~\ref{maximcl} that $X_2$ is closed.
\end{proof}  

\begin{defin} \label{topradical}
Let $X$ be a left unital Banach $A$-module. We say that the closed submodule of $X$ from
Theorem~\ref{exsrad} is the \emph{topological radical} of $X$ and denote it by $\trad X$. 
\end{defin}

\begin{prop}\label{irrtrad0}
The topological radical of an irreducible Banach module is trivial.
\end{prop}
\begin{proof}
 Let $X$ be an irreducible Banach module, and let $\phi$ be a small morphism with range in $X$. 
 Then $\Image\phi=X$ or $\Image\phi=0$. Since $\phi$ is small and  $X\not=0$,  $\phi$ is not surjective. Hence, $\phi=0$. This implies that $\trad X=0$.
\end{proof}   

\begin{prop}\label{fgrads}
Let $X$ be a finitely-generated Banach module. Then the natural embedding 
$\io\!:\trad X\to X$ is a small morphism. 
\end{prop}
\begin{proof}
Let $\phi$ be a  contractive monomorphism such that $\phi\dotplus\io$ is 
surjective. If $\phi$ is not surjective it follows from 
Proposition~\ref{fdmax} that there is a maximal  contractive monomorphism 
$\ga$ such that $\ga\succeq \phi$. Then $\ga\dotplus\io$ is surjective and 
$\ga\succeq \io$ (by the definition of the topological radical). This 
contradiction implies that $\phi$ is surjective. 
\end{proof} 

 \begin{prop}\label{fdmtr}
If $X$ is a unital finitely-generated Banach module over a unital Banach 
algebra $A$, then $\trad X=\rad X$. In particular, $\trad A$  coincides with 
the Jacobson radical of $A$.
\end{prop} 

The proposition follows immediately from Theorem~\ref{maximcl} and  the 
following lemma.

\begin{lem}\label{clmaxfd}
Every algebraically maximal submodule in a finitely-generated Banach module 
is closed.
\end{lem}
\begin{proof} 
Let $X_0$ be an algebraically maximal submodule in a finitely-generated 
Banach $A$-module $X$. Let $k$ be the minimal number such that for any 
finite set generating $X$ only $k$ generators are not contained in $X_0$. 
Note that $k>0$.

Fix generators $x_1,\ldots,x_n$ of $X$ such that $x_1,\ldots,x_k\in 
X\setminus X_0$ and $x_{k+1},\ldots,x_n\in X_0$. Denote by $U$ the set of 
all elements of the form $$x=\sum_{i>1}a_i\cdot x_i+(1-a_1)\cdot x_1,$$ 
where $\sum_{i}\|a_i\|<1$. 

If there exists $x_0\in X_0\cap U$, then $x_0=\sum_{i>1}a_i\cdot 
x_i+(1-a_1)\cdot x_1$, where $\|a_1\|<1$. Therefore, $1-a_1$ is invertible 
and 
$$
x_1=(1-a_1)^{-1}( x_0-\sum_{i>1}a_i\cdot x_i).
$$    
Hence, $x_0,x_2,\ldots,x_n$  are generators of $X$  but only $k-1$ 
generators are not in $X_0$. This contradiction with minimality of $k$ 
implies that $ X_0\cap U=\emptyset$.

It follows from the open mapping theorem that the surjective map $$A\ptn 
\ell^1_n\to X\!:e_i\mapsto x_i$$  is open. Therefore $U$ is open.  Since 
$x_1\in U$, we have $x_1 \not \in \overline{X_0}$.
  
Now assume that $\overline{X_0}\ne X_0$ and take $y \in 
\overline{X_0}\setminus X_0$. Since $X_0$ is maximal, we have $X_0+A\cdot 
y=X$. Therefore there are $a\in A$ and $x_0\in X_0$ such that $x_0+a\cdot 
y=x_1$. Note that $x_1+a\cdot(y'-y)=x_0+a\cdot y'\in X_0$ for every  $y'\in 
X_0$. However, since $x_1 \not \in \overline{X_0}$, we can take $y'\in X_0$ 
sufficiently close to $y$ to satisfy $x_1+a\cdot(y'-y) \not \in X_0$. This 
contradiction implies  that $X_0$ is closed.
\end{proof}

Now we can establish main properties of topological radical, which are 
similar to the algebraic case (cf. \cite[Sec.9.1, 9.2]{Ka}).

\begin{thm}\label{ptrstr} Let $X$ be a unital left Banach  $A$-module. 

 {\em (1)}  If $\phi\in \Ah(X,Y)$ then $\phi(\trad(X))\subset  \trad Y$. 
 
  {\em (2)}    $  \tau\!:A\to X\!:a\to a\cdot x_0$  is small iff  $x_0\in\trad X$.   
  
  {\em (3)}  $\overline{R\cdot X}\subset \trad X$, where $R$ is the Jacobson radical of $A$. 
  
  {\em (4)} $\trad(X/\trad X)=0$.  
  
  {\em (5)} If $Z$ is a closed submodule in $X$ such that $\trad(X/Z)=0$ then $\trad X\subset Z$. 
   \end{thm}     
\begin{proof}   
{\rm (1)} It follows from the definition and Proposition~\ref{compsm}.  
 
{\rm (2)} See the proof of  Theorem~\ref{exsrad}.
  
{\rm (3)}  Let  $x_0\in X$. It is sufficient to show that  $\tau'\!:R\to 
X\!:r\mapsto r\cdot x_0$ is small. Since $\tau'$ is a composition of  $R\to 
A$, which is small by Propositions~\ref{fgrads} and~\ref{fdmtr},  and 
$\tau\!:A\to X$,   Proposition~\ref{compsm} implies that $\tau'$ is small. 
  
{\rm (4)}  Suppose that $x\in X$ such that $x+\trad X\in  \trad(X/\trad X)$. 
By the definition $x+\trad X\in \Image \al$   for every maximal   
contractive monomorphism $\al\!:Y\to X/\trad X$. By Lemma~\ref{lifmax} there 
exists a maximal   contractive monomorphism $\be\!:W\to X$  such that $x\in 
\Image\be$. This implies
 that $x\in \trad X$. 
 
{\rm (5)}  Denote by $\si$ the projection $X\to X/Z$. It follows from (1) 
that $\si(\trad X)=0$. Therefore $\trad X\subset Z$.
\end{proof} 
   
\begin{cor}\label{rstr}
 $\rad X\subset\trad X$ for each unital left Banach  $A$-module $X$. 
   \end{cor}
 \begin{proof} 
Suppose that $x_0\in\rad X$. Then $A\cdot x_0$ is a small submodule in $X$ \cite[Sec.9.1.3(a)]{Ka}.  
Consider $\tau\!:A\to X\!:a\to a\cdot x_0$.
If $\phi\!:Y\to X$ such that $\phi\dotplus\tau$ is 
surjective then $X=A\cdot x_0+\Image\phi$. Since $A\cdot x_0$ is  small,
$\phi$ is surjective. Thus $\tau$ is a small morphism. By Theorem~\ref{ptrstr}(2)  $x_0\in\trad X$.
 \end{proof}

If  $A$ is not unital we can treat each Banach $A$-module $X$ as a unital 
Banach module over the unitization $A_+$ and consider the topological 
radical of $X$.

\begin{lem}\label{radbal}
Let $A$ be a radical Banach  algebra. Then 

{\em (1)}
$\rad A=A^2$ and  $\overline{A^2} \subset \trad A$;

{\em (2)} if $A$ admits a right b.a.i. $\rad A=A^2=\overline{A^2} =\trad A$.
\end{lem}
\begin{proof}
(1) Since $A_+/A$ is classically semi-simple, $A$ is a left good ring 
\cite[9.7.2,\,9.7.3(a)]{Ka}. Therefore $\rad X=A\cdot X$ for every left 
unital $A_+$- module $X$ \cite[9.7.1]{Ka}.  In particular, $\rad A=A^2$. The 
inclusion  $\overline{A^2} \subset \trad A$ follows  
Theorem~\ref{ptrstr}(3). 

The second statement follows from the Cohen factorization theorem.  
\end{proof}
Consider $C[0,1]$ and $L^1[0,1]$ as Banach algebras with respect to the 
cut-off convolution 
$$(f*g)(s)\!:=\int_0^s f(t)g(s-t)\,dt.$$ 
It is well known that both algebras are radical.

\begin{prop}\label{radco1}
{\em (1)}  If $A= (L^1[0,1],\ast)$, then $\trad A = \rad A$.

{\em (2)} If $A= (C[0,1],\ast)$, then $\trad A \ne \rad A$.
\end{prop}
\begin{proof}
(1)	Since $A=L^1[0,1]$ admits a b.a.i.,  Lemma~\ref{radbal} implies that $\trad A = \rad A$.

(2)	It is easy to see that $I_0\!:=\{f\in C[0,1]:\, f(0)=0\}$ is a closed ideal in $A= C[0,1]$ 
and $A^2\subset I_0$. Since smooth functions vanishing at $0$ are dense in $I_0$ and  
every such a function is a convolution of the derivative  and the constant, 
we have  $\overline{A^2} = I_0$. 
Note that $A/I_0$ is one-dimensional. 
This implies that $I_0\to A$ is a maximal contractive monomorphism. 
Therefore $\rad A\subset \trad A\subset I_0$. 
By Lemma~\ref{radbal}  $\rad A=A^2$ and  $I_0=\overline{A^2}=\trad A$. To see that $A^2\ne I_0$ note that
every function in $A^2$ is majorized by a linear function, therefore  $f(s)=\sqrt s$ is not in $A^2$.
\end{proof}

Recall that a left Banach $A$-module $P$ is called \emph{projective} if a 
morphism of Banach $A$-modules with range in $P$ admits a right inverse 
morphism provided it admits a right inverse bounded operator.

\begin{prop}\label{radpr}
If $P$ is a unital projective module with the approximation property, then 
$\trad P=\overline{R\cdot P}$, where $R$ is the Jacobson radical of~$A$.
\end{prop}
\begin{proof}
By Theorem~\ref{ptrstr}(3) $\overline{R\cdot P}\subset \trad P$. 

On the other hand, suppose that $x_0  \in \trad P$. Since $P$ is projective 
and has the approximation property,  \cite[Theorem~1(3)]{Se}  implies  that 
$x_0$ can be approximated in the norm topology by elements of the form 
$\sum_{i=1}^n\chi_i(x_0)\cdot y_i$ where $\chi_1,\ldots,\chi_n\in \Ah(P,A)$ 
and $y_1,\ldots,y_n\in P$. It follows from Theorem~\ref{ptrstr}(1) that 
 $\chi_i(x_0)\in R$.   Hence, $x_0\in \overline{R\cdot P}$.
\end{proof} 

\begin{re} 
It is not hard to check that in the case when $P$ is free, i.e. has the form 
$A\ptn E$ for some Banach space $E$, the argument of Proposition~\ref{radpr} 
can be applied to the case when $A$ or $E$ has the approximation property.
\end{re}   

\subsection*{Acknowledgements}
This research was partially  supported by RFBR under grant no.\ 12-01-00577.

\end{document}